\documentclass[11pt]{article}
\usepackage{amsthm,amssymb,amsmath,latexsym}
\usepackage{graphs}
\theoremstyle{definition}
\newtheorem{defin}{Definition}[section]
\theoremstyle{remark}

\newtheorem{remar}{Remark}
\theoremstyle{plain}
\newtheorem{thm}[defin]{Theorem}
\newtheorem{prop}[defin]{Proposition}
\newtheorem{lemm}[defin]{Lemma}
\newtheorem{corol}[defin]{Corollary}

\numberwithin{equation}{section}

\begin{document}

\title{Infinite measures on Cantor spaces}

\author{O.~Karpel\footnote{E-mail: helen.karpel@gmail.com }\\
\textit{Department of Mathematics}\\
\textit{Institute for Low Temperature Physics, Kharkov, Ukraine}}

\date{}
\maketitle

\begin{abstract}
We study the set $M_\infty(X)$ of all infinite full non-atomic Borel measures on a Cantor space $X$. For a measure $\mu$ from $M_\infty(X)$ we define a \textit{defective set} $\mathfrak{M}_\mu = \{x \in X : \mbox{for any clopen set } U \ni x \mbox{ we have } \mu(U) = \infty \}$. We call a measure $\mu$ from $M_\infty(X)$ \textit{non-defective} ($\mu \in M_\infty^0(X)$) if $\mu(\mathfrak{M}_\mu) = 0$.
The paper is devoted to the classification of measures $\mu$ from $M_\infty^0(X)$ with respect to a homeomorphism. The notions of goodness and clopen values set $S(\mu)$ are defined for a non-defective measure $\mu$. We give a criterion when two good non-defective measures are homeomorphic and prove that there exist continuum classes of weakly homeomorphic good non-defective measures on a Cantor space. For any group-like subset $D \subset [0,\infty)$ we find a good non-defective measure $\mu$ on a Cantor space $X$ with $S(\mu) = D$ and an aperiodic homeomorphism of $X$ which preserves $\mu$. The set $\mathcal{S}$ of infinite ergodic $\mathcal{R}$-invariant measures on non-simple stationary Bratteli diagrams consists of non-defective measures. For $\mu \in \mathcal{S}$ the set $S(\mu)$ is group-like, a criterion of goodness is proved for such measures. We show that a homeomorphism class of a good measure from $\mathcal{S}$ contains countably many distinct good measures from $\mathcal{S}$.


\end{abstract}

\section{Introduction}

In this paper, we initiate the classification of Borel infinite measures on a Cantor space with respect to a homeomorphism. Two measures $\mu$ and $\nu$ defined on Borel subsets of a topological space $X$ are called \textit{homeomorphic} if there exists a self-homeomorphism $h$ of  $X$  such that $\mu = \nu\circ h$, i.e. $\mu(E) = \nu(h(E))$ for every Borel subset $E$ of $X$. In such a way, the set of all Borel measures on $X$ is partitioned into equivalence classes.

The topological properties of the space $X$ are important for the classification of measures up to a homeomorphism. J.~Oxtoby and S.~Ulam \cite{Oxt-Ul} gave necessary and sufficient conditions under which a Borel probability measure on the finite-dimensional cube is homeomorphic to the Lebesgue measure. Later, similar results were obtained for various manifolds (see the book by S.~Alpern and V.~Prasad \cite{Alp-Pr} for the details).

A Cantor space is a non-empty zero-dimensional compact perfect metric space.
E.~Akin, T.~Austin, R.~Dougherty, R.~Mauldin, F.~Navarro-Bermudez, A.~Yingst (\cite{Akin3, Austin, D-M-Y, Navarro, Navarro-Oxtoby, Yingst}) studied the Borel \textit{probability} measures on Cantor spaces.
In those papers, the major results were focused on the classification of Bernoulli measures up to a homeomorphism.
It was E. Akin who started a systematic study of the Borel probability measures on a Cantor space~\cite{Akin1, Akin2}.
For Cantor spaces the situation is much more difficult than for connected spaces.
Though there is, up to a homeomorphism, only one Cantor space it is not hard to construct full (the measure of every non-empty open set is positive) non-atomic measures on the Cantor space which are not homeomorphic.
This fact is due to the existence of a countable base of clopen subsets of a Cantor space. For probability measures, the \textit{clopen values set} $S(\mu)$ is the set of values of a measure $\mu$ on all clopen subsets of $X$. This set provides an invariant for homeomorphic measures, although it is not a complete invariant, in general. But for the class of the so called \textit{good} probability measures, $S(\mu)$ \textit{is} a complete invariant. By definition, a  full non-atomic probability measure $\mu$ is good if whenever $U$, $V$ are clopen sets with $\mu(U) < \mu(V)$, there exists a clopen subset $W$ of $V$ such that $\mu(W) = \mu(U)$. It turns out that such measures are exactly invariant measures of uniquely ergodic minimal homeomorphisms of Cantor sets (see \cite{Akin2}, \cite{GW}).

Borel infinite measures arise as ergodic invariant measures for aperiodic homeomorphisms of a Cantor set. The study of homeomorphic infinite measures is of crucial importance for the classification of Cantor aperiodic systems up to orbit equivalence.
Let $\mu$ be an infinite Borel measure on a Cantor space $X$. Denote by $\mathfrak{M}_\mu$ the set of all points in $X$ whose clopen neighbourhoods have only infinite measures. Then $\mathfrak{M}_\mu$ is a closed non-empty subset of $X$.
We study the full non-atomic infinite measures $\mu$ such that $\mu(\mathfrak{M}_\mu) = 0$. We call such measures \textit{non-defective}.
A full ergodic infinite Borel measure $\mu$ is either non-defective or it gives every non-empty clopen set infinite measure.

Section 2 of the work is focused on the study of non-defective measures on Cantor spaces in general while Section 3 is focused on the particular class of non-defective measures which are ergodic invariant measures on non-simple Bratteli diagrams.

In Section 2, we prove that a non-defective measure $\mu$ is $\sigma$-finite and $X = \bigsqcup_{i = 1}^\infty V_i \bigsqcup \mathfrak{M}_\mu$, where each $V_i$ is a clopen set of finite measure and $\mathfrak{M}_\mu$ is closed, nowhere dense and has zero measure. Thus, non-defective measures are uniquely determined by their values on clopen subsets of $X$.
Though for a non-defective measure $\mu$ there always exist countably many clopen sets of infinite measure, we define the clopen values set $S(\mu)$ to be the set of all finite values of the measure $\mu$ on the clopen sets. We define a good non-defective measure in the same way as a good probability measure. We show that some theorems of~\cite{Akin1, Akin2, Akin3} for probability measures have analogues for non-defective measures. Besides there are specific results concerning infinite measures.

It turns out that in the case of non-defective measures, clopen values set is not a complete invariant for homeomorphisms of good measures. Let $\mu$ and $\nu$ be good non-defective measures on Cantor spaces $X$ and $Y$. Let $\mathfrak{M}$ be the defective set for $\mu$ and $\mathfrak{N}$ be the defective set for $\nu$. Then $\mu$ is homeomorphic to $\nu$ if and only if $S(\mu) = S(\nu)$ and there exists a homeomorphism $\widetilde{h} \colon \mathfrak{M} \rightarrow \mathfrak{N}$ where the sets $\mathfrak{M}$ and $\mathfrak{N}$ are endowed with the induced topologies.
For a good probability measure $\mu_1$ the group of self-homeomorphisms of $X$ that preserve $\mu_1$ acts transitively on $X$. We show that the group of self-homeomorphisms of $X$ that preserve a good non-defective measure $\mu$ acts transitively on $X \setminus \mathfrak{M}_\mu$ and it cannot act transitively on $X$.
A subset $D$ of $[0,\infty)$ is \textit{group-like} when $D = G \bigcap [0,\infty)$ with $G$ an additive subgroup of reals. We show that for every countable dense group-like subset $D$ of $[0,\infty)$ there are good non-homeomorphic non-defective measures $\mu$ and $\nu$ on a Cantor space such that $S(\mu) = S(\nu) = D$. In the probability case for any group-like subset $D_1$ of the unit interval there exists a unique up to homeomorphism good measure $\mu_1$ with $S(\mu_1) = D_1$.
We call two non-defective measures $\mu$ and $\nu$ on a Cantor space $X$ \textit{weakly homeomorphic} if there exists a self-homeomorphism $h$ of $X$ and a constant $C>0$ such that $\mu \circ h = C \nu$. We prove that there exist continuum classes of good weakly homeomorphic non-defective measures on Cantor spaces.

We provide several examples of non-defective measures on Cantor spaces. In Section 2, it is shown how to obtain a non-defective measure by extending an infinite locally finite measure from a locally compact space $Y$ to its compactification $cY$.
In Section 3, we consider non-simple stationary Bratteli diagrams and non-simple Bratteli diagrams of finite rank. The infinite ergodic measures on their path spaces invariant with respect to the cofinal (tail) equivalence relation $\mathcal{R}$ are studied. These measures are also non-defective.
The article~\cite{S.B.} is used, where an explicit description of all ergodic (finite and infinite) measures on stationary Bratteli diagrams is found. We also use~\cite{S.B.O.K.} where homeomorphic probability measures on stationary Bratteli diagrams are studied.
We show that the clopen values set is a complete invariant for homeomorphism equivalence of infinite good ergodic $\mathcal{R}$-invariant measures on stationary Bratteli diagrams.
For any infinite ergodic invariant measure $\mu$ on a stationary Bratteli diagram it is proved that the set $S(\mu)$ is group-like. We give a criterion of goodness for $\mu$.
It is proved that for every good ergodic $\mathcal{R}$-invariant measure $\mu\in \mathcal S$ there exist countably many good ergodic $R$-invariant measures $\{\mu_i\}_{i\in \mathbb N}$ on stationary Bratteli diagrams such that the measures $\mu$ and $\mu_i$ are homeomorphic but the tail equivalence relations on the corresponding Bratteli diagrams are not orbit equivalent.
We also use~\cite{BKMS} where it is shown how infinite measures on Bratteli diagrams of finite rank arise from finite measures. We show how to build a good measure on a Bratteli diagram of finite rank with any given countable dense group-like subset $D$ of $\mathbb{Q}$.

\section{Infinite measures on Cantor spaces}\label{AnSurv}
In this section, we study a class of non-atomic infinite Borel measures on a Cantor space.
We define the notion of good infinite measure and prove the analogues of the theorems for good probability measures from~\cite{Akin2, Akin3}.

Let $X$ be a Cantor space.
A measure on a Cantor space is called \textit{full} if every non-empty open set has positive measure. Recall that a \textit{support} of a measure $\mu$ is a closed set $supp(\mu)$ such that $\mu(X \setminus supp(\mu)) = 0$ and $\mu(supp(\mu) \cap U) > 0$ for every open $U$ such that $supp(\mu) \cap U$ is non-empty. It is easy to see that for a non-atomic measure $\mu$ the set $supp(\mu)$ in the induced topology is a Cantor set. We can consider measures on their supports to obtain full measures.
Denote by $M_\infty(X)$ the set of all full non-atomic Borel infinite measures on $X$.

\begin{defin}
Let $\mu \in M_\infty(X)$. Define $\mathfrak{M}_\mu = \{x \in X : \mbox{for any clopen set } U \ni x \mbox{ we have } \mu(U) = \infty \}$. We call the set $\mathfrak{M}_\mu$ \textit{defective}.
\end{defin}

We will need the following properties of the set $\mathfrak{M}_{\mu}$.

\begin{prop}\label{setM}
Let $\mu \in M_\infty(X)$. Then the set $\mathfrak{M}_{\mu}$ is a non-empty closed subset of $X$.
\end{prop}

\noindent \textbf{Proof.} Assume that $\mathfrak{M}_{\mu} = \emptyset$. Then every point $x \in X$ has a clopen neighbourhood $U_x$ such that $\mu(U_x) < \infty$. Since $X$ is compact, it follows that $\mu(X) < \infty$. This is a contradiction.

To prove that $\mathfrak{M}_\mu$ is closed, consider a set $X \setminus \mathfrak{M}_\mu = \{x \in X : \mbox{ there exists a clopen set } U_x \ni x \mbox{ such that } \mu(U_x) < \infty \}$. Then for every point $x \in X \setminus \mathfrak{M}_\mu$ we have $U_x \subset X \setminus \mathfrak{M}_\mu$. Hence $X \setminus \mathfrak{M}_\mu$ is open.
$\blacksquare$

\begin{defin}
Denote by $M^{0}_\infty(X)$ the class of measures $\mu \in M_\infty(X)$ such that  $\mu(\mathfrak{M}_\mu) = 0$. We call a measure from the set $M^{0}_\infty(X)$ \textit{non-defective}.
\end{defin}

\begin{prop}\label{strangeproperty}
Let $\mu \in M^{0}_\infty(X)$. Then

(1) For any clopen set $U$ with $\mu(U) = \infty$ and any $a > 0$ there exists a clopen subset $V \subset U$ such that $a \leq \mu(V) < \infty$.

(2) $X = \bigsqcup_{i = 1}^\infty V_i \bigsqcup \mathfrak{M}_\mu$, where each $V_i$ is a clopen set of finite measure and $\mathfrak{M}_\mu$ is closed nowhere dense and has zero measure.

(3) $\mu$ is $\sigma$-finite

(4) $\mu$ is uniquely determined by its values on the algebra of clopen sets.
\end{prop}

\noindent \textbf{Proof.} (1) Let $U$ be a non-empty clopen subset of $X$ such that $\mu(U) = \infty$. Since $\mu \in M^{0}_\infty(X)$, we have $\mu(U) = \mu(U \setminus \mathfrak{M}_{\mu})$.
There are only countably many clopen subsets in $X$, hence the open set $U \setminus \mathfrak{M}_{\mu}$ can be represented as a disjoint union of clopen subsets $\{U_i\}_{i \in \mathbb{N}}$ of finite measure.
We have $\mu(U) = \sum_{i=0}^{\infty} \mu(U_i) = \infty$, hence for every $a \in \mathbb{R}$ there is a clopen subset $V = \bigsqcup_{i=0}^{N}U_i$ such that $a \leq \mu(V) < \infty$.

(2) Suppose that $\mathfrak{M}_\mu$ is dense in some clopen set $U$. Then $U \subset \mathfrak{M}_\mu$ since $\mathfrak{M}_\mu$ is closed. On the one hand, $\mu(U) \leq \mu(\mathfrak{M}_\mu) = 0$. On the other hand, $U$ is a clopen neighbourhood of its points, hence $\mu(U) = \infty$ by definition of $\mathfrak{M}_\mu$.

(3) and (4) Obviously follow from (2).
$\blacksquare$

\begin{remar}
Our aim is to study ergodic infinite measures invariant with respect to aperiodic homeomorphisms of a Cantor set. Let $\mu \in M_\infty(X)$. Let $h$ be a self-homeomorphism of $X$ which preserves measure $\mu$. Then $h(\mathfrak{M}_\mu) = \mathfrak{M}_\mu$. If $\mu$ is ergodic for $h$ then either $\mu(\mathfrak{M}_\mu) = 0$ (hence $\mu \in M_\infty^{0}(X)$) or $\mu(X \setminus \mathfrak{M}_\mu) = 0$.
The equality $\mu(X \setminus \mathfrak{M}_\mu) = 0$ implies that
every clopen subset of $X$ has zero or infinite measure.
\end{remar}

Section 3 of the paper is devoted to the study of ergodic invariant measures on stationary Bratteli diagrams and Bratteli diagrams of finite rank. Every such measure $\mu$ on its support $supp(\mu)$ belongs to $M^{0}_\infty(supp(\mu))$ (see~\cite{S.B., BKMS}). Another example of non-defective measures is given in the proof of Theorem~\ref{goodSmu}. There the non-defective measures are obtained from finite measures on a Cantor space $X$ by compactification of the space $X \times \mathbb{Z}$.

Recall that a $\sigma$-additive measure $\nu$ on $\mathfrak{B}$, or a measure space $(X,\mathfrak{B},\nu)$, is called \textit{semi-finite} if whenever $E \in \mathfrak{B}$ and $\nu(E) = \infty$ there is an $F \subseteq E$ such that $F \in \mathfrak{B}$ and $0 < \nu(F) < \infty$. If a measure is $\sigma$-finite then it is semi-finite (see~\cite{Fr}).

The following proposition provides other examples of the measures from $M_\infty^0(X)$.
\begin{prop}\label{extension}
Let $\mathfrak{C}$ be the algebra of all clopen subsets of $X$. Let $\mu$ be a finitely additive measure on $\mathfrak{C}$ such that

(1) $\mu$ is non-atomic, i.e. for every $U \in \mathfrak{C}$ such that $\mu(U) > 0$ there exists $V \in \mathfrak{C}$ such that $V \subset U$ and $0 < \mu(V) < \mu(U)$,

(2) $\mu(U) > 0$ for every non-empty $U \in \mathfrak{C}$,

(3) $\mu(X) = \infty$,

(4) If $\mu(U) = \infty$ for $U \in \mathfrak{C}$ and $a \in \mathbb{R}$ then there is $V \subset U$ such that $V \in \mathfrak{C}$ and $a \leq \mu(V) < \infty$.

Then $\mu$ has a unique extension to a Borel measure from $M_\infty^0(X)$.
\end{prop}

\noindent \textbf{Proof.} We use the well-known construction of outer measure. Set $\widehat{\mu}(A) = \inf \{\sum_{j = 0}^{\infty} \mu(C_i) : C_i \in \mathfrak{C} \mbox{ and } A \subseteq \bigcup C_i\}$ for any $A \subset X$. Since $\mathfrak{B}$ is generated by $\mathfrak{C}$, the outer measure $\widehat{\mu}$ is a $\sigma$-additive measure on $\mathfrak{B}$.
The measure space $(X, \mathfrak{B}, \widehat{\mu})$ is not semi-finite. Indeed, every point in $\mathfrak{M}_{\widehat{\mu}}$ has infinite measure. In such cases we can use the semi-finite version of $(X,\mathfrak{B}, \widehat{\mu})$ (see~\cite{Fr}). For $E \in \mathfrak{B}$ set $\mu_{sf}(E) = \sup\{\widehat{\mu}(E \bigcap F) : F \in \mathfrak{B},\; \widehat{\mu}(F) < \infty\}$. Then $(X,\mathfrak{B}, \mu_{sf})$ is a semi-finite measure space. The values of $\mu_{sf}$ on the clopen sets are the same as of $\mu$. It is not hard to prove that $\mu_{sf} \in M_\infty^0(X)$. $\blacksquare$

Since the set $\mathfrak{M}_\mu$ is non-empty, there exist countably many clopen sets of infinite measures. For a measure $\mu \in M^0_\infty(X)$ define the \textit{clopen values set} as the set of all finite values of the measure $\mu$ on the clopen sets:
$$
S(\mu) = \{\mu(U):\,U\mbox{ is clopen in } X \mbox{ and } \mu(U) < \infty\}.
$$

For each measure $\mu \in M_\infty^0(X)$, the set $S(\mu)$ is a countable dense subset of the interval $[0, \infty)$. Indeed, by Proposition~\ref{strangeproperty}, for every $a \in [0,\infty)$ there is a clopen subset $V$ such that $a \leq \mu(V) < \infty$. The set $S(\mu)$ is dense in $[0, \mu(V)]$ (see~\cite{Akin1}).

Let $X_{1}$, $X_{2}$ be two Cantor sets.
It is said that measures $\mu_{1}$ on $X_{1}$ and $\mu_{2}$ on $X_{2}$ are \textit{homeomorphic} if there exists a homeomorphism $h \colon X_{1} \rightarrow X_{2}$ such that $h_{*}\mu_{1} = \mu_{2}$. Clearly, $S(\mu_{1}) = S(\mu_{2})$ for any homeomorphic measures $\mu_1$ and $\mu_2$.
We call two Borel infinite measures $\mu$ and $\nu$ on a Cantor space $X$ \textit{weakly homeomorphic} if there exists a self-homeomorphism $h$ of $X$ and a constant $C>0$ such that $h_{*} \mu = C \nu$.

\medskip
Let $D$ be a dense countable unbounded subset of the interval $[0,\infty) \subset \mathbb{R}$.
Then $D$ is called \textit{group-like} if for every $\gamma \in D$ there exists an additive subgroup $G$ of $\mathbb{R}$ such that $D \bigcap \; [0,\gamma] = G \bigcap \; [0,\gamma]$. It is clear that for all positive $\gamma \in D$ the group $G$ is the same.
For $G$ an additive subgroup of $\mathbb{R}$ we call a positive real number $\alpha$ a \textit{divisor} of $G$ if $\alpha G = G$. The set of all divisors of $G$ is called $Div(G)$.
For any $A \subset \mathbb{R}$ let $G[A]$ denote the additive group generated by A, so that $G[A]$ is the set of all finite sums of differences of elements of $A$ (including $0 =$ the empty sum so that $G[\emptyset] = \{0\}$).

\begin{lemm}\label{grouplike_subsets}
Let $D$ be a dense countable unbounded subset of the interval $[0, \infty)$ with $0 \in D$.
The following conditions on $D$ are equivalent:

(1) $D$ is group-like.

(2) $D = G \bigcap [0,\infty)$ for $G$ an additive subgroup of $\mathbb{R}$.

(3) $(D \bigcap [0, \gamma]) + \gamma \mathbb{Z}$ is an additive subgroup of $\mathbb{R}$ for every $\gamma \in D$.

(4) $\alpha, \beta \in D$ and $\alpha \leq \beta$ imply that $\beta - \alpha \in D$.

\end{lemm}

\noindent \textbf{Proof.} Implications $(2)\Rightarrow (1)$ and $(1)\Rightarrow (4)$ are clear.

$(4)\Rightarrow (3)$: Following~\cite{Akin2} we show that $(D \bigcap [0, \gamma]) + \gamma \mathbb{Z}$ is closed under subtraction. Let $\alpha, \beta \in (D \bigcap [0, \gamma]) + \gamma \mathbb{Z}$. Then $\alpha = d_1 + \gamma n_1$ and $\beta = d_2 + \gamma n_2$ for some $d_1, d_2 \in D \bigcap [0, \gamma]$ and $n_1, n_2 \in \mathbb{Z}$. Then $\alpha - \beta = d_1 - d_2 + \gamma (n_1 - n_2)$. If $d_1 - d_2 > 0$ then $d_1 - d_2 \in D \bigcap [0, \gamma]$ by (4). Otherwise $d_1 - d_2 + \gamma = \gamma - (d_2 - d_1) \in D \bigcap [0, \gamma]$ and $\alpha - \beta = d_2 - d_1 + \gamma (n_1 - n_2 + 1)$. Hence, in any case, $\alpha - \beta \in (D \bigcap [0, \gamma]) + \gamma \mathbb{Z}$.

$(3)\Rightarrow (2)$: Let $\gamma \in D$. We note that if $(D \bigcap [0, \gamma]) + \gamma \mathbb{Z}$ is an additive group for every $\gamma \in D$ then it equals $G[D]$. Indeed, suppose for some $\gamma_1, \gamma_2 \in D$ with $\gamma_1 < \gamma_2$ we have $(D \bigcap [0, \gamma_1]) + \gamma_1 \mathbb{Z} = G_1$ and $(D \bigcap [0, \gamma_2]) + \gamma_2 \mathbb{Z} = G_2$. Then $G_1 \bigcap [0, \gamma_1] = G_2 \bigcap [0, \gamma_1]$ and $G_1 = G_2 = G[D]$.
Thus, for every $\gamma \in D$, $(D \bigcap [0, \gamma]) + \gamma \mathbb{Z} = G$, where $G$ is an additive subgroup of $\mathbb{R}$. We have $G \bigcap [0, \infty) = \bigcup_{\gamma \in D} D \bigcap [0, \gamma] = D$. $\blacksquare$

\medskip
We give below the  definitions of \textit{good, refinable, and weakly refinable measures} which are based on the corresponding definitions for probability measures.
We follow here the papers~\cite{Akin2, Akin3,D-M-Y}.

A \textit{partition basis} $\mathcal{B}$ for a Cantor set $X$ is a collection of clopen subsets of $X$ such that every non-empty clopen subset of $X$ can be partitioned by elements of $\mathcal{B}$. A partition basis is a basis for the topology but not every basis is a partition basis (see~\cite{Akin3}).

\begin{defin}
Let $\mu$ belong to $M_\infty^0(X)$.

(1) A clopen subset $V$ of $X$ is called \textit{good} for $\mu$ (or just good when the measure is understood) if for every clopen subset $U$ of $X$ with $\mu(U) < \mu(V)$, there exists a clopen set $W$ such that $W \subset V$ and $\mu(W) = \mu(U)$. A measure $\mu$ is called \textit{good} if every clopen subset of $X$ is good for $\mu$.

(2) A clopen subset $U$ of $X$ of finite measure is called \textit{refinable} for $\mu$ if $\alpha_{1}, \ldots ,\alpha_{k} \in S(\mu)$ with $\alpha_{1}+ \ldots + \alpha_{k} = \mu(U)$ implies that there exists a clopen partition $\{U_{1}, \ldots , U_{k}\}$ of $U$ with $\mu(U_{i}) = \alpha_{i}$ for $i = 1,\ldots ,k$.
We call a measure $\mu$ \textit{refinable} if every clopen subset of finite measure is refinable and every clopen subset of infinite measure is good.

(3) A measure $\mu \in M_\infty(X)$ is called \textit{weakly refinable} if there exists a partition basis $\mathcal{B}$ for $X$ with $X \in \mathcal{B}$ consisting of refinable and good clopen subsets.
\end{defin}

\begin{remar}
A measure $\mu \in M^0_\infty(X)$ is good if every clopen subset of finite measure is good. Indeed, since $\mu$ belongs to $M^0_\infty(X)$, the clopen set of infinite measure contains good clopen subsets of arbitrary large measures.
As in the case of probability measures, the notions of refinability and weakly refinability are the weakenings of the notion of goodness. We need to perform algebraic operations with $\infty$, such as $\infty + \infty = \infty$, to define refinable clopen set of infinite measure similarly to the refinable clopen set of finite measure. In Corollary~\ref{2mes} there is an example of a good measure $\mu \in M_\infty^0(X)$ such that $\mathfrak{M}_\mu$ consists of one point. Hence any partition of a clopen set $V$ with $\mu(V) = \infty$ contains only one subset of infinite measure.
\end{remar}

The proofs of the following proposition and corollary are similar to the proofs of Proposition 2.11 and Corollary 2.12 in~\cite{Akin3}.
\begin{prop}
Let $\mu \in M_\infty^0(X)$.

(a) If $S(\mu)$ is group-like and a clopen $U \subset X$ of finite measure is refinable then $U$ is good.

(b) If the measure $\mu$ is good then every clopen subset of $X$ with finite measure is refinable.

\end{prop}

\begin{corol}\label{goodmeasure}
For a measure $\mu \in M^0_\infty(X)$ the following are equivalent:

(a) $\mu$ is a good measure.

(b) $\mu$ is refinable and $S(\mu)$ is group-like.

(c) $\mu$ is weakly refinable and $S(\mu)$ is group-like.
\end{corol}

\begin{thm}\label{Partition_Basis}
Let $\mu$ belong to $M^0_\infty(X)$. Then $\mu$ is good if and only if there exists a partition basis $\mathcal B$ consisting of clopen sets which are good for $\mu$. In particular, if a clopen set can be partitioned by good clopen sets, then it is itself good.
\end{thm}

\noindent \textbf{Proof.} The ``only if'' part of the theorem is trivial.
We prove the ``if'' part. Given two clopen sets $U,V$ with $\mu(U) < \mu(V)$ we have to construct a clopen subset $W$ of $V$ such that $\mu(W) = \mu(U)$. If $\mu(V) < \infty$, the proof is contained in~\cite{Akin3}, see Theorem 2.7. If $\mu(V) = \infty$ then there exists a partition $\{B_1,...,B_k\}$ of $V$ into good clopen sets. Evidently, there exists $B_i$ among $\{B_1,...,B_k\}$ such that $\mu(B_i) = \infty$. Since $B_i$ is good, there exists $W \subset B_i$ with $\mu(W) = \mu(U)$. \quad $\blacksquare$

\begin{thm}\label{product}
If $\mu \in M^0_\infty(X)$, $\nu \in M^0_\infty(Y)$ are good measures, then the product $\mu \times \nu$ is a good measure on $X \times Y$ and
$$
S(\mu \times \nu) = \left\{\sum_{i=0}^N \alpha_i \cdot \beta_i : \alpha_i \in S(\mu), \beta_i \in S(\nu), N \in \mathbb{N}\right\}.
$$
\end{thm}

\noindent \textbf{Proof.} Since $\mu\times \nu (\mathfrak{M_{\mu\times \nu}}) \leq \mu\times \nu(\mathfrak{M}_\mu \times Y) + \mu\times \nu(X \times \mathfrak{M}_\nu) = 0$, we see that $\mu \times \nu \in M^0_\infty(X \times Y)$.
If $U \subset X$ and $V \subset Y$ are clopens then $U \times V$ is a clopen subset of $X \times Y$. We call such set a rectangular clopen. In~\cite{Akin3}, it is shown that rectangular clopens form a partition basis for $X \times Y$. By Theorem~\ref{Partition_Basis}, it suffices to show that an arbitrary non-empty clopen $U \times V$ is good for $\mu \times \nu$. Given finite $\alpha_1,...,\alpha_k \in S(\mu)$, $\beta_1,...,\beta_k \in S(\nu)$ such that $\sum_i \alpha_i \cdot \beta_i < \mu(U) \cdot \nu(V)$, we must prove that there exists a clopen subset of $U \times V$ of measure $\sum_i  \alpha_i \cdot \beta_i$.
If $\mu(U) < \infty$ and $\nu(V) < \infty$, then we use the proof of Theorem 2.8 in~\cite{Akin3}. If $\mu(U) = \infty$, $\nu(V) < \infty$ then, by Proposition~\ref{strangeproperty}, there exists a clopen subset $W \subset U$ of finite measure such that $\sum_i \alpha_i \cdot \beta_i < \mu(W) \cdot \nu(V)$. The same idea works in the case when both measures are infinite.
$\blacksquare$

The following corollary can be proved in the same way as Corollary 2.9 in~\cite{Akin3}:
\begin{corol}
The product of finite or infinite sequence of good measures is a good measure.
\end{corol}

The following theorem is one of the main results of the work. It provides the criterion when two non-defective measures are homeomorphic. The analogous theorem was proved in~\cite{Akin3} for probability measures.

\begin{thm}\label{weakref_homeo}
Let $\mu$ and $\nu$ be good non-defective measures on Cantor spaces $X$ and $Y$. Let $\mathfrak{M}$ be the defective set for $\mu$ and $\mathfrak{N}$ be the defective set for $\nu$. Then $\mu$ is homeomorphic to $\nu$ if and only if the following conditions hold:

(1) $S(\mu) = S(\nu)$,

(2) There exists a homeomorphism $\widetilde{h} \colon \mathfrak{M} \rightarrow \mathfrak{N}$ where the sets $\mathfrak{M}$ and $\mathfrak{N}$ are endowed with the induced topologies.

\end{thm}

\noindent \textbf{Proof}.
The ``only if'' part is trivial.

To prove the ``if'' part, we use back and forth construction described in Theorem 2.14 of~\cite{Akin3}. Denote the clopen subsets of $X$ as $\{C_0, C_1,...\}$ and the clopen subsets of $Y$ as $\{D_0, D_1,...\}$. We begin with one-element partitions $\mathcal{X}_0 = \{X\}$ and $\mathcal{Y}_0 = \{Y\}$, and the map $\rho_0$ which sends element $\mathcal{X}_0$ to $\mathcal{Y}_0$.

Suppose we have defined partitions $\mathcal{X}_{2k}$ and $\mathcal{Y}_{2k}$ of $X$ and $Y$ respectively and a bijection $\rho_{2k} \colon \mathcal{X}_{2k} \rightarrow \mathcal{Y}_{2k}$ so that

(i)\; $\nu(\rho_{2k}(A)) = \mu(A) \mbox{ for all } A \in \mathcal{X}_{2k}$,

(ii)\; $\rho_{2k}(A) \bigcap \mathfrak{N} = \widetilde{h}(A \bigcap \mathfrak{M})$.

Let $\mathcal{X}_{2k+1} = \{A \bigcap C_k : A \in \mathcal{X}_{2k}\} \bigcup \{A \setminus C_k : A \in \mathcal{X}_{2k}\}$. Then $\mathcal{X}_{2k+1}$ is a partition of $X$ into clopen sets and $\mathcal{X}_{2k+1}$ refines $\mathcal{X}_{2k}$. Fix some $B \in \mathcal{Y}_{2k}$. Those elements of $\mathcal{X}_{2k+1}$ which are subsets of $\rho_{2k}^{-1}(B)$ form a partition of $\rho_{2k}^{-1}(B)$. Denote this partition by $\{A_1,...,A_n\}$. Let $\mu(A_i) < \infty$ for $i = 1,...,m$ and $\mu(A_i) = \infty$ for $i = m + 1,...,n$.

Since $\nu$ is good and $S(\mu) = S(\nu)$, we can choose disjoint clopen subsets $B_1,..,B_m \subset B$ such that $\nu(B_i) = \mu(A_i)$, $i = 1,...,m$.
If $m = n$ then $B_1,..,B_m \subset B$ form a partition of $B$. Suppose $m < n$. By (ii), we have $\widetilde{h}(A_j \bigcap \mathfrak{M}) \subset B \setminus \bigsqcup_{i=1}^m B_i$ for $j = m+1,....,n$. Since $\widetilde{h}(A_j \bigcap \mathfrak{M})$ is a clopen subset of $\mathfrak{N}$ in the induced topology, there exists a clopen set $U_j \subset Y$ such that $U_j \bigcap \mathfrak{N} = \widetilde{h}(A_j \bigcap \mathfrak{M})$. Set $V_j = U_j \bigcap (B \setminus \bigsqcup_{i=1}^m B_i)$. Then $V_j \subset B \setminus \bigsqcup_{i=1}^m B_i$ is a clopen subset of infinite measure and
\begin{equation}\label{inf}
\widetilde{h}(A_j \bigcap \mathfrak{M}) = V_j \bigcap \mathfrak{N}.
\end{equation}

Since $\{A_j\}_{j = m+1}^n$ are disjoint, we see that the sets $\{V_j \bigcap \mathfrak{N}\}_{j = m+1}^n$ are disjoint. Therefore, $\mu(V_j \bigcap V_k) < \infty$. Moreover, $\mu(B \setminus \bigcup_{j=m+1}^n V_j) < \infty$. Hence we can make the sets $\{V_j\}_{j = m+1}^n$ be a disjoint partition of $B \setminus \bigsqcup_{i=1}^m B_i$ and preserve the equality~(\ref{inf}). Indeed, set $B_{m+1} = V_{m+1}$ and $B_j = V_j \setminus \bigcup_{k = m+1}^{j-1} V_k$ for $j = m+1,...,n-1$. Set $B_n = B \setminus \bigsqcup_{i = 1}^{n-1}B_i$. Then the sets $\{B_i\}_{i=1}^n$ form a partition of $B$ such that $\nu(B_i) = \mu(A_i)$ and $\widetilde{h}(A_j \bigcap \mathfrak{M}) = B_j \bigcap \mathfrak{N}$ for $i = 1,...,n$.

Let $\mathcal{Y}_{2k+1}$ include $\{B_1,...,B_n\}$, and let $\rho_{2k+1}$ map $A_i$ to $B_i$ for $i = 1..n$. After doing this for each $B$ in $\mathcal{Y}_{2k}$, we obtain $\mathcal{Y}_{2k+1}$, a refinement of $\mathcal{Y}_{2k}$, and a bijection $\rho_{2k+1}$ such that conditions (i), (ii) hold.

Now reverse the roles of $X$ and $Y$. We may repeat the above construction to define $\mathcal{X}_{2k+2}$, $\mathcal{Y}_{2k+2}$ and $\rho_{2k+2}$.

We have defined $\mathcal{X}_{k}$, $\mathcal{Y}_{k}$ and $\rho_k$ for all $k \geq 0$. For each $x \in X$ and for each $k \geq 0$ there is a unique $A_k \in \mathcal{X}_{k}$ such that $x \in A_k$, since $\mathcal{X}_{k}$ is a partition of $X$. By construction, the sequence $\{A_k\}_{k\geq 0}$ will be a nested sequence of clopen sets. Let $h$ be a map from $X$ to $Y$ such that $h(x)$ is an element of $\bigcap_{k=0}^{\infty} \rho_k(A_k)$, another nested sequence. This element is unique since the sets in $Y_{2k+2}$ separate points in $D_k$ from points not in $D_k$, and every clopen set is one of these $D_k$'s. It is easy to verify that the map $h$ is a homeomorphism.
Since we consider Borel infinite measures which are uniquely defined by their values on the clopen sets, $h$ preserves the measure of every Borel set.
$\blacksquare$

\begin{remar}
Theorem~\ref{weakref_homeo} is also true for weakly refinable measures. Recall that by definition the clopen sets of infinite measure are good for weakly refinable measure. The proof of the theorem changes only in the cases where clopen sets of finite measures are considered. In those cases we use the proof of Theorem  2.14 in~\cite{Akin3}.
\end{remar}

\begin{remar}
The analogous theorem can be formulated and proved for weakly homeomorphic measures:

\textit{Let $\mu$ and $\nu$ be good non-defective measures on Cantor spaces $X$ and $Y$. Let $\mathfrak{M}$ be the defective set for $\mu$ and $\mathfrak{N}$ be the defective set for $\nu$. Then $\mu$ is weakly homeomorphic to $\nu$ if and only if the following conditions hold:}

\textit{(1) There exists $c > 0$ such that $S(\mu) = c S(\nu)$,}

\textit{(2) There exists a homeomorphism $\widetilde{h} \colon \mathfrak{M} \rightarrow \mathfrak{N}$ where the sets $\mathfrak{M}$ and $\mathfrak{N}$ are endowed with the induced topologies.}

\end{remar}

\begin{remar}
Let $\mu \in M^0_{\infty}(X)$ be a good measure and $V$ be any clopen subset of $X$ with $\mu(V) < \infty$. Then $\mu$ on $X$ is homeomorphic to $\mu$ on $X \setminus V$.
Let $S(\mu) = G \bigcap [0, \infty)$. Then $\mu$ is homeomorphic to $c \mu$ if and only if $c \in Div(G)$.

\end{remar}

Denote by $H_{\mu}(X)$ the group of measure preserving homeomorphisms of a Cantor space $X$. The action of $H_{\mu}(X)$ on $X$ is called \textit{transitive} if for every $x_{1}, x_{2} \in X$ there exists $h \in H_{\mu}(X)$ such that $h(x_{1}) =  x_{2}$. The action is called \textit{topologically transitive} if there exists a dense orbit, i.e. there is $x \in X$ such that the set $O(x) =  \{h(x) : h \in H_\mu(X)\}$ is dense in $X$.

\begin{corol}
Let $\mu$ be a good measure in $M^0_{\infty}(X)$. Then the group $H_\mu(X)$ acts transitively on $X \setminus \mathfrak{M_\mu}$. In particular, the group $H_\mu(X)$ acts topologically transitively on $X$.
\end{corol}

\noindent \textbf{Proof.} Let $x_1, x_2 \in X\setminus \mathfrak{M_\mu}$. Since $\mu$ is good, there exist disjoint clopen sets $U_1$, $U_2$ such that $\mu(U_1) = \mu(U_2) < \infty$ and $x_i \in U_i$ for $i = 1,2$. By Proposition 2.11 of~\cite{Akin2}, there exists a homeomorphism $h_1 \colon U_1 \rightarrow U_2$ such that $h(x_1) = x_2$ and $h$ preserves $\mu$. By Theorem~\ref{weakref_homeo}, there exists a homeomorphism $h_2 \colon X \setminus U_1 \rightarrow X \setminus U_2$ which preserves $\mu$. Define $h \colon X \rightarrow X$ to be $h_1$ on $U_1$ and $h_2$ on $X \setminus U_1$. Clearly, $h$ is a measure-preserving homeomorphism such that $h(x_1) = x_2$.

For every $x \in X\setminus\mathfrak{M}$ we have $O(x) = X \setminus\mathfrak{M}_\mu$. By Proposition~\ref{setM}, the set $X \setminus\mathfrak{M}_\mu$ is dense in $X$. Hence $H_\mu(X)$ acts topologically transitively on $X$. $\blacksquare$

\begin{thm}\label{goodSmu}
Let $\mu \in M_\infty^{0}(X)$ be a good measure. Then the clopen values set $S(\mu)$ is group-like. Conversely, for every countable dense group-like subset $D$ of $[0, \infty)$ which contains $0$ there is a good measure $\mu$ on Cantor space such that $S(\mu) = D$.
\end{thm}

\noindent \textbf{Proof.} The first statement of the theorem can be proved the same way as Proposition 2.4 in~\cite{Akin2}. If $D$ is group-like and $\gamma \in D$ then $\frac{1}{\gamma}D \bigcap [0,1]$ is a group-like subset of $[0,1]$. In~\cite{Akin2} it was proved that there exists a good probability measure $\mu_1$ on Cantor space $X$ with $S(\mu_1) = \frac{1}{\gamma}D \bigcap [0,1]$. Set $\mu = \gamma \mu_1$. Then $\mu$ is good and $S(\mu) = D \bigcap [0,\gamma]$.

Endow the set $\mathbb{Z}$ with discrete topology. Let $\alpha\mathbb{Z} = \mathbb{Z} \bigcup \{\infty\}$ be a one-point compactification of $\mathbb{Z}$. Let $\nu$ be the extension of the counting measure on $\mathbb{Z}$ such that $\nu(\{\infty\})= 0$.
Consider product measure $\widetilde{\mu} = \mu \times \nu$ on a Cantor space $\widetilde{X} = X \times \alpha\mathbb{Z}$. The measure $\widetilde{\mu}$ is an infinite full nonatomic Borel measure on $\widetilde{X}$. It is clear that $\mathfrak{M}_{\widetilde{\mu}} = X \times \{\infty\}$ and $\widetilde{\mu}(\mathfrak{M}_{\widetilde{\mu}}) = 0$.

Since the measure $\mu$ is good it easily follows that the measure $\widetilde{\mu}$ is good. Indeed, it suffices to prove that every clopen subset of finite measure is good. The rectangular clopen sets $\{U \times \{z\}\}$, where $U$ is clopen in $X$ and $z \in \mathbb{Z}$ form the partition basis for clopen sets of finite measure. Since $\widetilde{\mu} (U \times \{z\}) = \mu(U)$, these sets are good.

We have $S(\widetilde{\mu}) = \{\sum_{i=1}^n \mu(U_i) : U_i \mbox{ is clopen in } X, \; i = 1,...,n \}$. Since $D$ is group-like, we see that $S(\widetilde{\mu}) = D$.
$\blacksquare$

\medskip
The following corollary shows that unlike the case of probability measures, there are non-homeomorphic good non-defective measures with the same clopen values set.
\begin{corol}\label{2mes}
For every countable dense group-like subset $D$ of $[0, \infty)$ there exist two non-homeomorphic good non-defective measures $\mu_1$ and $\mu_2$ on Cantor spaces $X_1$ and $X_2$ such that $S(\mu_1)= S(\mu_2) = D$.
\end{corol}
\noindent \textbf{Proof.}
Let $X$ be a Cantor space.
Let $X_1 = X \times \alpha\mathbb{Z}$ and $\mu_1$ be as in the proof of Theorem~\ref{goodSmu}. Every non-compact locally compact Hausdorff topological space has a one-point compactification. Let $X_2$ be a one-point compactification of the space $X\times\mathbb{Z}$, i.e. $X_2 = \alpha(X \times \mathbb{Z}) = (X\times\mathbb{Z}) \bigcup \{\infty\}$. Let $\mu_2 = \mu \times \nu$ on $X\times\mathbb{Z}$ and $\mu_2(\{\infty\}) = 0$. Then $\mu_1$ and $\mu_2$ are good measures and $S(\mu_1) = S(\mu_2) = D$. The set $\mathfrak{M}_{\mu_1}$ is a Cantor set while the set $\mathfrak{M}_{\mu_2}$ consists of one point $\{\infty\}$. Hence, by Theorem~\ref{weakref_homeo}, the measures $\mu_1$ and $\mu_2$ are not homeomorphic. $\blacksquare$

\medskip
Note that an infinite Borel measure cannot be a unique invariant measure for a homeomorphism of a Cantor space, since every homeomorphism of a compact metric space has a non-trivial set of invariant Borel probability measures.

\begin{corol}\label{invarhomeo}
Let $D$ be a countable dense group-like subset of $[0, \infty)$. Then there exists an aperiodic homeomorphism of a Cantor space with good non-defective invariant measure $\widetilde{\mu}$ such that $S(\widetilde{\mu}) = D$.
\end{corol}

\noindent \textbf{Proof}. We use the construction from the proof of Theorem~\ref{goodSmu}. For given set $D$ we build a good infinite measure $\widetilde{\mu} = \mu \times \nu$ on $\widetilde{X} = X \times \alpha\mathbb{Z}$ such that $S(\widetilde{\mu}) = D$. Recall that $\mu$ is a good measure on the Cantor space $X$ with $S(\mu) = D \bigcap [0,\gamma]$ for some $\gamma \in D$ and $\nu$ is an extension of a counting measure on $\mathbb{Z}$. Since the measure $\mu$ is a good finite measure on $X$, there exists a minimal homeomorphism $T$ of $X$ such that $\mu$ is invariant for $T$ (see~\cite{Akin2}). Let $T_1(x,n) = (Tx, n+1)$. Then $T_1$ is aperiodic homeomorphism of $\widetilde{X}$. The measure $\widetilde{\mu}$ is invariant for $T_1$. $\blacksquare$

\begin{thm}
There exist continuum distinct classes of weakly homeomorphic good measures in $M_{\infty}^{0}(X)$.
\end{thm}

\noindent \textbf{Proof.} Let $\{g_\alpha\}_{\alpha \in \Lambda}$ be a Hamel basis for $\mathbb{R}$ over $\mathbb{Q}$. Then it contains one rational element which we denote by $g_\lambda$.
Consider additive subgroups $G_\alpha = G[1, g_\alpha]$ of $\mathbb{R}$ for $\alpha \in \Lambda \setminus\{\lambda\}$. 
Then each $G_\alpha$ is a countable dense subgroup of $\mathbb{R}$. Obviously we have $G_\alpha = G_\beta$ if and only if $\alpha = \beta$.

Consider an equivalence relation on the set $\{G_\alpha\}_{\alpha \in \Lambda \setminus \{\lambda\}}$. We call two groups $G_\alpha$ and $G_\beta$ equivalent if there exists $c > 0$ such that $G_\alpha = c \; G_\beta$. Then the set $\{G_\alpha\}_{\alpha \in \Lambda\setminus \{\lambda\}}$ is partitioned into equivalence classes. We show that every equivalence class contains at most countable number of elements. Therefore, there are continuum such equivalence classes.

Suppose that $G_\alpha = c \;G_\beta$. Then we have $c \in G_\alpha$ and $\beta \in \frac{1}{c} G_\alpha$.
Since $G_\alpha$ is countable, there are at most countably many numbers $c$ and $\beta$ such that the equality holds. Hence for every $G_\alpha$ there exist at most countably many groups from  $\{G_\beta\}_{\beta \in \Lambda\setminus \{\lambda\}}$ such that $G_\alpha = c \; G_\beta$ for some $c \in \mathbb{R}$.
By Theorem~\ref{goodSmu}, for every $G_\alpha$ there exists a good measure $\mu_\alpha \in M_\infty^{0}(X)$ such that $S(\mu) = G_\alpha \bigcap [0,\infty)$. By Theorem~\ref{weakref_homeo}, there exist continuum distinct classes of weakly homeomorphic measures. $\blacksquare$

\medskip
If a Cantor space $X$ is ordered, then the \textit{special clopen values set} $\widetilde{S}(\mu)$ can be defined. For probability measures, the set $\widetilde{S}(\mu)$ is a complete invariant with respect to a measure preserving and order preserving homeomorphism(see~\cite{Akin1, Akin2}).
We show that for infinite measures the notion of order is much harder to use.

By an \textit{order} $\leq$ on a space $X$ we mean a total order, i.e. any two points are comparable and the order is closed as a subset of $X \times X$. Equivalently, the order topology is the original compact topology on $X$. An ordered Cantor space is a pair $(X, \leq)$ where $X$ is a Cantor space and $\leq$ is an order on $X$.
For an ordered Cantor space $(X, \leq)$ there exists a minimal element which we will denote by $m$ and a maximal element which we will denote by $M$.
We will adopt the interval notation so that, for example, $[x,y) = \{z \in X : x \leq z < y\}$. A point $x \in X$ is called a \textit{left endpoint} if the interval $[m,x]$ is clopen and a \textit{right endpoint} if $[x, M]$ is clopen.
If $(X, \leq)$ is an ordered space, then the \textit{special clopen values set} for $\mu$ is
$$
\widetilde{S}(\mu) = \{\mu([m,x]) : x \mbox{ is an endpoint of X}, \mu([m,x]) < \infty \}.
$$

Clearly, $\widetilde{S}(\mu) \subset S(\mu)$.

\begin{defin}
Let $(X, \leq)$ be an ordered Cantor space and $\mu \in M^0_{\infty}(X)$. Let $m$ be the minimal element of $X$.

1. We call the order $\leq$ \textit{good} if for any two points $x,y \in X$ if $\mu([x,y]) < \infty$ then $\mu([m,x]) < \infty$ and $\mu([m,y]) < \infty$.

2. The measure $\mu$ is called \textit{adapted} to $(X, \leq)$, an order $\leq$ is called adapted to $\mu$, when
$$
\widetilde{S}(\mu) = S(\mu).
$$
\end{defin}

\begin{remar}
For a good order we have $\mu([x,y]) = \mu([m,y]) - \mu([m,x])$ for any set $[x,y]$ of finite measure. Hence we have $S(\mu)\subset G[\widetilde{S}(\mu)]$.
\end{remar}

Since for a non-defective measure $\mu$ the set $\mathfrak{M}_\mu$ is nowhere dense, it is not hard to prove the following theorem, we leave the proof to a reader.
\begin{thm}\label{onepointM}
Let $\mu \in M^0_\infty(X)$. Then there exists a good order for $\mu$ if and only if the set $\mathfrak{M}_\mu$ consists of one point.
\end{thm}

The following theorem can be proved the same way as Theorem 2.6 in~\cite{Akin2}.
\begin{thm}

Let $(X, \leq)$ be an ordered Cantor space and $\mu \in M^0_\infty(X)$.

(a) If the measure $\mu$ is adapted to $(X, \leq)$ then the set $\widetilde{S}(\mu)$ is a group-like subset of $[0, \infty)$. The converse is true for a good order $\leq$.

(b) If $\mu$ is adapted to $(X, \leq)$, then there exists a clopen subset $U$ with $\mu(U) = \infty$ such that $U$ is good for $\mu$. Moreover, if the order $\leq$ is good then the measure $\mu$ is good.

\end{thm}

\begin{thm} \label{eq_implies_good}
Let $\mu$ be a non-defective measure on a Cantor space $X$.
If for every $x$ from $X \setminus\mathfrak{M_\mu}$ there exists an order $\leq$ such that $x$ is a minimal point of $(X, \leq)$ and $S(\mu) = \widetilde{S}(\mu)$, then $X$ has a basis of good clopen sets.
\end{thm}

\noindent \textbf{Proof.} Since we consider the order which induces topology, every non-empty clopen set of $X$ can be written as a finite disjoint union of clopen intervals (see~\cite{Akin1}):
$$
V = \bigsqcup_{i=1}^k V_i \mbox{ with } V_i = [x_i, y_i],
$$

Let $z$ be any point in $X \setminus \mathfrak{M}_\mu$. Then there exists an order $\leq$ such that $z$ is its minimum point and $\widetilde{S}(\mu) = S(\mu)$. Then any clopen set of the kind $[z,x]$ is good.

Hence for any clopen set $V$ and any point $x \in V \setminus \mathfrak{M}_\mu$ there exists an order $\leq_x$ and a good clopen subset $[x,y] \subset V$. So any point $x \in V \setminus \mathfrak{M}_\mu$ has a good clopen neighbourhood. The set $X \setminus \mathfrak{M}_\mu$ is dense in $X$. Hence the set $V$ can be covered by good clopen sets.
$\blacksquare$

\begin{thm}\label{good_implies_equality}
Let $\mu$ be a good non-defective measure on a Cantor space $X$. Let $x_0$ belong to $X \setminus \mathfrak{M}_\mu$. Then there exists an order $\leq$ on $X$ such that the minimal point $m$ of $(X, \leq)$ is $x_0$ and the measure $\mu$ is adapted to $(X, \leq)$, i.e. $\widetilde{S}(\mu) = S(\mu)$.
\end{thm}

\noindent \textbf{Proof.} The proof is analogous to the proof of Theorem 2.7 in~\cite{Akin2}. The only difference is that we should change the definition of $mesh_\mu(\mathcal{A})$ for a partition $\mathcal{A}$ of the space $X$. Namely, we set $mesh_\mu(\mathcal{A}) = max\{\mu(A): A \in \mathcal{A} \mbox{ and } \mu(A) < \infty\}$.
$\blacksquare$

\begin{remar}
In the case of probability measures the measure $\mu_1$ on $X$ is good if and only if there exists an order $\leq$ on $X$ such that $\widetilde{S}(\mu) = S(\mu)$.
For a good probability measure $\nu$ on the Cantor space $X$ and for any point $x_0 \in X$ there exists an order $\leq$ such that $\nu$ is adapted to $(X,\leq)$ and $x_0$ is a minimal point (see~\cite{Akin2}). In the case of infinite measure $\mu$ there always exists a point $y \in X$ such that for any order $\leq$ with the minimal point $y$ we have $\widetilde{S}(\mu) = \{0\}$. Such point $y$ should be taken from the set $\mathfrak{M}_\mu$.
\end{remar}

\section{Infinite Measures on Bratteli Diagrams}\label{section3}

We will use the notation and results from~\cite{S.B.O.K., S.B.}.
For basic definitions and facts about Bratteli diagrams see \cite{S.B.,GPS}.
Let $B$ be a stationary non-simple Bratteli diagram
and $A = F^T$ be the matrix transposed to incidence matrix of the diagram.
In~\cite{S.B.}, it was shown that any ergodic (both finite and infinite) measure is an extension of a finite ergodic measure from a simple subdiagram. Recall that a real number $\lambda$ is called \textit{distinguished} if there exists a non-negative eigenvector $x$ with $Fx = \lambda x$. If $\lambda$ is a distinguished Perron-Frobenius eigenvalue for $A$ then the corresponding $\mathcal R$-invariant measure on $X_B$ is finite, otherwise it is infinite~(see \cite{S.B.}). We consider measures on $X_B$ only on their supporting sets.

Let $\mu$ be a full infinite ergodic $\mathcal{R}$-invariant measure on a stationary Bratteli diagram $B$ defined by a non-distinguished eigenvalue $\lambda$ of the matrix $A = F^T$. Let $\alpha$ be the class of vertices corresponding to $\mu$.
Let the vertex set of subdiagram $B_f$ consist of vertices of the class $\alpha$ and of vertices of the classes $\beta \prec \alpha$ such that the measures of cylinder sets that end in the vertices of class $\beta$ are finite. Let the edge set of $B_f$ consist of all edges of $B$ that start and end in the vertices of $B_f$. Denote by $A_f = F_f^T$ the matrix transpose to the incidence matrix of $B_f$. Then $A_f x = \lambda x$.
Denote by $\mu_f$ the probability ergodic invariant measure on $B_f$ generated by $x$, $\lambda$. Then $\lambda$ is a distinguished Perron-Frobenius eigenvalue for $A_f$ and $x$ is a reduced vector for $\mu_f$ (see~\cite{S.B.O.K.}). Denote by $X_{B_f}$ the set of all infinite paths of $B_f$. Clearly $X_{B_f}$ supports $\mu_f$.

\begin{remar}
We have introduced a sort of ``normalizing'' for full infinite ergodic $\mathcal{R}$-invariant measures on a stationary Bratteli diagram. Indeed, let $\nu$ be a full infinite ergodic $\mathcal{R}$-invariant measure on a stationary Bratteli diagram $B$. By normalizing $\nu_f$ to be a probability measure, we choose an element of the class of weakly homeomorphic measures for a measure $\nu$.
\end{remar}

Recall that $X_\alpha$ is the set of all infinite paths in $X_B$ that eventually pass only through the vertices of the class $\alpha$.
We have $\mathfrak{M}_\mu = X_B \setminus X_\alpha$ and $\mathfrak{M}_\mu$ is a Cantor space in the induced topology. Since $\mu$ is ergodic we have $\mu(\mathfrak{M}_\mu) = 0$ and $\mu \in M_\infty^0(X_B)$. Since all Cantor spaces are homeomorphic, by Theorem~\ref{weakref_homeo}, the clopen values set is a complete invariant for good ergodic $\mathfrak{R}$-invariant measures on stationary Bratteli diagrams.

Let $\gamma$ be the vertex class of the diagram such that $\gamma \prec \alpha$ and the measures of cylinder sets that end in the vertices of class $\gamma$ are infinite. Then we call $\gamma$ an ``infinite'' component of $B$. All infinite paths of $X_B$ that eventually pass only through the vertices of the class $\gamma$ belong to $\mathfrak{M}_\mu$.
Denote by $\widetilde{h}^{(N)}_i$ be the number of the cylinder sets that end in the vertex $i$ on the level $N$ and do not belong to $X_{B_f}$.

\begin{thm}\label{grouplikeinf}
Let $\mu$ be an infinite ergodic $\mathcal{R}$-invariant measure on a stationary diagram $B$ defined by a non-distinguished eigenvalue $\lambda$ of the matrix $A = F^T$. Let $(x_1, \ldots, x_n)^T$ be the corresponding reduced vector and $H$ the additive subgroup of $\mathbb{R}$ generated by $\{x_1, \ldots , x_n\}$. Then the clopen values set $S(\mu)$ is group-like and
$$
S(\mu) = \left(\bigcup_{N=0}^\infty \frac{1}{\lambda^N} H\right) \cap [0,\infty).
$$
\end{thm}

\noindent\textbf{Proof.}
We prove that for any $m \in H(x_1,...,x_n) \bigcap [0,\infty)$ and any $N \in \mathbb{N}$ there exists a clopen set $U \subset X_B$ such that
$$
\mu(U) = \frac{m}{\lambda^N}.
$$
Let the "infinite" vertex class of the diagram be accessed from the class $\beta \preceq \alpha$. Let the vertices with measures $x_{l+1},...,x_{l+r}$ belong to the class $\beta$.
Denote by
$$
L = \left\{\sum_{i = l+1}^{l+r}k_i^{(M)}\frac{x_i}{\lambda^{M-1}} : 0 \leq k_i^{(M)} \leq \widetilde{h}^{(M)}_i; \;M = 1,2,...\right\}.
$$
Then $S(\mu) = S(\mu_f) + L$.
For any element $s \in L$ there exists a clopen set $W \subset X_B \backslash X_B^f$ such that $\mu(W) = s$. Since the measure $\mu$ is infinite, the set $L$ contains arbitrary large elements.

Since $0 < \frac{x_i}{\lambda^M} < 1$ for $i = 1,...,n$ and $M \in \mathbb{N}$, there exists $s \in L$ such that $\frac{m}{\lambda^N} - s \in (0,1)$. But the element $\frac{m}{\lambda^N} - s$ belongs to $\left(\bigcup_{N=0}^\infty \frac{1}{\lambda^N} H\right) \bigcap [0,1]$. Hence it belongs to $S(\mu_f)$. Therefore there exists a clopen set $V \subset X_B^f$ such that $\mu(V) = \frac{m}{\lambda^N} - s$. Set $U = V \bigsqcup W$. Then $U$ is the clopen subset of $X_B$ and $\mu(U) = \frac{m}{\lambda^N}$.
$\blacksquare$

From Corollary~\ref{goodmeasure} and Theorem~\ref{grouplikeinf} we obtain the following
\begin{corol}
For an ergodic invariant measure $\mu$ on a stationary Bratteli diagram the following are equivalent:

(i) $\mu$ is a good measure.

(ii) $\mu$ is refinable.

(iii) $\mu$ is weakly refinable.

\end{corol}

\begin{lemm}\label{gdinf}
Let $B$ be a stationary Bratteli diagram. Let $\mu$ be an ergodic infinite $\mathcal R$-invariant measure on $B$. Denote by $\alpha$ the non-distinguished class of vertices that defines $\mu$. Then $\mu$ is good if and only if all the clopen cylinder sets that end in the vertices of the class $\alpha$ are good.
\end{lemm}

\noindent\textbf{Proof.}
The proof is similar to the proof of Lemma 3.4 in~\cite{S.B.O.K.}.
$\blacksquare$

\begin{thm}\label{KritGoodinf}
Let $B$ be a stationary Bratteli diagram. Let $\mu$ be an ergodic infinite $\mathcal R$-invariant measure on $B$. Let $\mu_f$ be the corresponding finite measure. Then $\mu$ is good if and only if $\mu_f$ is good.
\end{thm}

\noindent\textbf{Proof.}
To prove the ``if'' part we use Lemma~\ref{gdinf} and the ideas of the proof of Theorem 3.5 in~\cite{S.B.O.K.}. Let $\mu_{fin}$ be good. Suppose the cylinder set $U$ ends in the vertex of the class $\alpha$. If $U$ belongs to $X_{fin}$ then it is good since $\mu_{fin}$ is good. Else if $U$ belongs to $X_{inf}$ then there exists a cylinder set $V \subset X_{fin}$ such that $V$ is tail equivalent to $U$ (and hence $\mu(V) = \mu(U)$). The set $V$ is good, hence for every clopen set $W$ with $\mu(W) < \mu(V)$ there exists a clopen subset $Y \subset V$ such that $\mu(Y) = \mu(W)$. Since $V$ is tail equivalent to $U$ there exists a clopen subset $Z \subset U$ such that $Z$ is tail equivalent to $Y$. Hence for every cylinder set $U$ ends in the vertex of the class $\alpha$ and any clopen set $W$ with $\mu(W) < \mu(U)$ there exists a clopen subset $Z \subset U$ such that $\mu(Z) = \mu(W)$. By Lemma~\ref{gdinf}, the measure $\mu$ is good.

The ``only if'' part of the result is obvious.
$\blacksquare$

\begin{remar}
It is easy ti see that Lemma~\ref{gdinf} and Theorem~\ref{KritGoodinf} are true for Bratteli diagrams of finite rank (see~\cite{S.B.O.K.,BKMS}).
\end{remar}

\begin{corol}\label{kriGood}
Let $\mu$ be an infinite $\mathcal R$-invariant measure on a stationary diagram $B$ defined by a non-distinguished eigenvalue $\lambda$ of the matrix $A = F^T$. Denote by $x = (x_1,...,x_n)^T$ the corresponding reduced vector. Let the vertices $m+1, \ldots, n$ belong to the class $\alpha$ corresponding to $\mu$. Then $\mu$ is good if and only if there exists $R \in \mathbb{N}$ such that $\lambda^R x_1,...,\lambda^R x_m$ belong to the additive group generated by $\{x_j\}_{j=m+1}^n$.

If $S(\mu)$ consists of rational numbers and $(\frac{p_1}{q}, \ldots, \frac{p_n}{q})^T$ is the corresponding reduced vector, then $\mu$ is good if and only if $\gcd(p_{m+1}, ...,p_n) | \; \lambda^R$ for some $R \in \mathbb N$.
\end{corol}

Let $D \subset [0, \infty)$ be an unbounded set. Then a positive integer $n$ is called a \textit{reciprocal} for $D$ if $\frac{1}{n} \in D$. Let $Rec(D)$ denote the set of reciprocals of $D$.

In the next theorem we show that for every countable dense unbounded group-like subset $D$ of $[0, \infty)$ there is a good infinite ergodic $\mathcal{R}$-invariant measure $\widehat{\mu}$ on the path space of a Bratteli diagram of finite rank such that $S(\widehat{\mu}) = D$. We will need the following modification of Lemma 2.2 of Akin~\cite{Akin2}, the proof is analogous:

\begin{lemm}\label{rec}
Let $D$ be a group-like subset of $[0, \infty)$ with $0,1 \in D$.

If $\frac{m}{n} \in D$ with $m,n$ relatively prime positive integers then $n \in Rec(D)$. If $n \in Rec(D)$ and $m | n$, then $m \in Rec(D)$. If $m, n \in Rec(D)$, then the least common multiple, $lcm(m,n)$, is in $Rec(D)$.

In general,
$$
\mathbb{Q} \bigcap D = \left\{\frac{k}{n} : k \in \mathbb{N} \mbox{ and } n \in Rec(D) \right\}.
$$
If $Rec(D)$ is infinite, then $\mathbb{Q} \bigcap D$ is dense in $[0, \infty)$. If $Rec(D)$ is finite, then $\mathbb{Q} \bigcap D \bigcap [0,1]$ is finite.
\end{lemm}

Now we can prove the following theorem.
\begin{thm}
Let $D$ be a dense group-like subset of $\mathbb{Q} \bigcap [0, \infty)$. Then there exists a good infinite ergodic $\mathcal{R}$-invariant measure $\widehat{\mu}$ on a Bratteli diagram of finite rank such that $S(\widehat{\mu}) = D$.
\end{thm}

\noindent \textbf{Proof}. First we construct a probability measure $\mu$ on the odometer such that $S(\mu) = D \bigcap [0, 1]$. Recall that an odometer can be represented as a Bratteli diagram $B$ which has one vertex on each level. Enumerate the reciprocals of $D$ such that $Rec(D) = \{r_i\}_{i = 0}^{\infty}$, $r_i = r_j$ if and only if $i = j$. Let $P = \{p_j\}_{j = 0}^{\infty}$ be the set of prime divisors of $\{r_i\}_{i = 0}^{\infty}$ such that every divisor $p_j$ is written as many times as $t(p_i) = \max_{r \in Rec(D)}\{\mbox{ multiplicity of } p_j \mbox{ in } r\}$. The number $t(p_i)$ can be either finite or infinite.

Let the number of the edges between neighbouring vertices $v_i$ and $v_{i+1}$ of the diagram $B$ be $p_i$. Consider the ergodic probability measure $\mu$ on $X_B$. The measure $\mu$ is good since it is the only invariant probability measure for a minimal homeomorphism of a Cantor space~\cite{Akin2}. It is clear that
$$
S(\mu) = \left\{\frac{k}{p_{1} p_2 \ldots p_{N}}, k, N \in \mathbb{N},\; 0 \leq k \leq p_{1} \ldots p_{N}\right\}.
$$
By Lemma~\ref{rec}, $S(\mu) = D \bigcap [0,1]$. Indeed, every denominator of the elements of $S(\mu)$ can be written as $d = \prod_{i = 1}^m q_i^{\alpha_i}$. Here $q_i \in P$ and $q_i \neq q_j$ for $i \neq j$. It follows from the construction that for every $i$ there exists $n_i \in Rec(D)$ such that $q_i^{\alpha_i} | n_i$. Hence, by Lemma~\ref{rec}, every $q_i^{\alpha_i} \in Rec(D)$. Since the numbers $\{q_i\}_{i = 1}^m$ are relatively prime, we obtain that $lcm(\{q_i^{\alpha_i}\}_{i = 1}^m) = \prod_{i = 1}^m q_i^{\alpha_i}$. Hence $d \in Rec(D)$ and $S(\mu) \subset D \bigcap [0,1]$. On the other hand, any element from $Rec(D)$ can be written as $\prod_{i = 1}^m p_i^{\alpha_i}$, where $\alpha_i \leq t(p_i)$. Hence $S(\mu) = D \bigcap [0,1]$.

Now we construct the non-simple Bratteli diagram $\widehat{B}$ with an infinite ergodic $\mathcal{R}$-invariant measure $\widehat{\mu}$ such that $\widehat{\mu}$ is an extension of $\mu$ and $S(\widehat{\mu}) = S(\mu) + \mathbb{N} = D$.

Let $\widehat{B}$ be defined by the sequence of matrices
$$
F_N =
\begin{pmatrix}
p_n & 0\\
a_n & p_n \\
\end{pmatrix}.
$$
Thus, $\widehat{B}$ has two components: $\beta \prec \alpha$, which look exactly like diagram $B$.
In~\cite{BKMS} there are given sufficient conditions for a measure $\widehat{\mu}$ to be infinite.
The numbers $a_n$ must satisfy the following equality:
$$
\sum_{i = 1}^{\infty} \frac{a_n}{p_n} = \infty.
$$
For instance we can choose $a_n = p_n$. The measure of every cylinder set that ends in the vertices of the class $\beta$ is infinite. Every clopen set is a finite union of cylinder sets. Hence any clopen set $U$ of finite measure is a finite union of cylinders that end in the vertices of class $\alpha$. Therefore, $\widehat{\mu}(U) \in G[S(\mu)]$. The fact that $\widehat{\mu}$ is good and $S(\widehat{\mu})$ is group-like can be proved as in Theorems~\ref{grouplikeinf},~\ref{KritGoodinf}.
$\blacksquare$

\begin{thm}\label{goodhomeoinf}
Let $\mu$ be a good ergodic $\mathcal{R}$-invariant infinite measure on a stationary (non-simple) Bratteli diagram $B$. Then there exist stationary Bratteli diagrams  $\{B_i\}_{i=0}^\infty$ and  good ergodic $\mathcal{R}_i$-invariant probability measures $\mu_i$ on $B_i$ such that each measure $\mu_i$ is homeomorphic to $\mu$ and the dynamical systems $(B_i, \mathcal R_i)$, $(B_j, \mathcal R_j)$ are topologically orbit equivalent if and only if $i = j$. Moreover, the diagram $B_i$ has exactly $i$ minimal components for the tail equivalence relation $\mathcal R_i, i\in \mathbb N$.
\end{thm}

\noindent \textbf{Proof.} We can add arbitrary many ``infinite'' minimal components to the given Bratteli diagram. This operation preserves $S(\mu)$ and the goodness of the measure. $\blacksquare$

\begin{prop}
Let $\mu$ be a good ergodic $\mathcal{R}$-invariant infinite measure on a stationary Bratteli diagram $B$. Then $\mu$ is homeomorphic to the measure $\nu$ for which there exists an aperiodic homeomorphism $T$ with exactly one minimal component and $\nu$ is the only ergodic infinite invariant measure for $T$.
\end{prop}

\noindent \textbf{Proof.} In Theorem 4.1 in~\cite{S.B.O.K.}, it is proved that there exists a stationary simple Bratteli diagram $B_{fin}$ with probability invariant measure $\mu_{fin}$ such that $S(\mu_{fin}) = S(\mu) \bigcap [0,1]$. By adding an "infinite" minimal component to $B$ we obtain a stationary Bratteli diagram $B'$ with the only infinite ergodic measure $\nu$. From Theorems~\ref{grouplikeinf} and~\ref{KritGoodinf}, it follows that $\nu$ is good and $S(\mu) = S(\nu)$. By Theorem~\ref{weakref_homeo}, the measures $\mu$ and $\nu$ are homeomorphic.
$\blacksquare$

\begin{prop}
For any ergodic $\mathcal{R}$-invariant measure on a stationary Bratteli diagram, the good order does not exist.
\end{prop}

\noindent \textbf{Proof.} The proposition follows from Theorem~\ref{onepointM}. $\blacksquare$

\medskip
Recall that a partition basis is a basis for the topology but not every basis is a partition basis.
\begin{thm}\label{goodonBD}
In order that an ergodic $\mathcal{R}$-invariant measure $\mu$ (either probability or infinite) on a stationary Bratteli diagram be good, it suffices that there exists a basis $\mathcal B$ consisting of clopen sets which are good for $\mu$. In particular, if a clopen set can be covered by good clopen sets, then it is itself good.
\end{thm}

\noindent \textbf{Proof.} Assume the converse. Then there exists a basis $\mathcal B$ consisting of clopen sets which are good for $\mu$, but the measure $\mu$ is bad. Hence every clopen subset of $\mu$ contains a good clopen subset. Let $\mu$ be defined by a non-distinguished eigenvalue $\lambda$ of the matrix $A = F^T$. Denote by $x = (x_1,...,x_n)^T$ the corresponding reduced vector. Let the vertices $m+1, \ldots, n$ belong to the class $\alpha$ corresponding to $\mu$. By Corollary~\ref{kriGood}, for any $R \in \mathbb{N}$ there exists $i \in \{1,...,m\}$ such that $\lambda^R x_i \not \in H(x_{m+1},...,x_n)$. Since $A (x_{m+1},...,x_n)^T = \lambda (x_{m+1},...,x_n)^T$, we get $\lambda H(x_{m+1},...,x_n) \subset H(x_{m+1},...,x_n)$. Therefore, if $\lambda^{R_0} x_i \in H(x_{m+1},...,x_n)$ for some $R_0 \in \mathbb{N}$, then $\lambda^R x_i \in H(x_{m+1},...,x_n)$ for all $R > R_0$. Thus we obtain that there exists $i \in \{1,...,m\}$ such that for any $R \in \mathbb{N}$ we have $\lambda^R x_i \not \in H(x_{m+1},...,x_n)$.

Consider any cylinder set $V$ that ends in the vertex of the class $\alpha$. We obviously have $\mu(V) = \frac{x_j}{\lambda^N}$ for some $j \in \{m+1,...,n\}$ and $N \in \mathbb{N}$. Hence $\mu(V) \in  H(\frac{x_{m+1}}{\lambda^N},...,\frac{x_n}{\lambda^N})$. Any clopen subset of $V$ can be represented as a finite disjoint union of cylinders that end in the vertices of the class $\alpha$. Hence the measure of any clopen subset of $V$ lies in $\bigcup_{K \in \mathbb{N}} H(\frac{x_{m+1}}{\lambda^K},...,\frac{x_n}{\lambda^K})$.

There exists $M \in \mathbb{N}$ such that $\frac{x_i}{\lambda^M} < \frac{x_j}{\lambda^N}$. Recall that $\frac{x_i}{\lambda^M}$ does not belong to $H(\frac{x_{m+1}}{\lambda^K},...,\frac{x_n}{\lambda^K})$ for any $K \in \mathbb{N}$. So $\frac{x_i}{\lambda^M} \not \in \bigcup_{K \in \mathbb{N}} H(\frac{x_{m+1}}{\lambda^K},...,\frac{x_n}{\lambda^K})$.
Besides, the number $\frac{x_i}{\lambda^M}$ is a measure of a cylinder set. Hence the set $V$ is bad and similarly all its clopen subsets are bad. We get a contradiction. $\blacksquare$

From Theorems~\ref{eq_implies_good}, \ref{goodonBD} we obtain the following corollary:

\begin{corol}
An ergodic infinite $\mathcal{R}$-invariant measure $\mu$ on a stationary Bratteli diagram $B$ is good if and only if for every $x \in X_B \setminus \mathfrak{M_\mu}$ there exists an order $\leq$ on $X_B$ such that $\widetilde{S}(\mu) = S(\mu)$.
\end{corol}

We give several examples illustrating the results of this section.

\noindent \textbf{Example 1.} Consider a class of stationary Bratteli diagrams $B_N$ with incidence matrices
$$
F_N =
\begin{pmatrix}
M(N) & 0 & 0 & 0\\
1 & 2 & 0 & 0\\
0 & 1 & N & 1\\
0 & 1 & 1 & N \\
\end{pmatrix},
$$
where $N \geq 3$ and $M(N) > N$. Let $\mu_N$ be a full infinite measure on $B_N$. It is generated by the non-distinguished eigenvalue $\lambda_N = N+1$. Let $\mu_{fin}(N)$ be the maximal finite measure for $\mu_N$. In~\cite{S.B.O.K.} it was shown that $\mu_{fin}(N)$ is a good measure if and only if for all sufficiently large $R$ we have $\frac{(N+1)^R}{N} \in \frac{N-1}{2N}\mathbb Z$. Thus, $\mu_N$ is good under the same conditions. For instance, $\mu_N$ is good for $N = 3$ but is not good for $N = 4$. This gives us an example of two infinite measures $\mu_3$ and $\mu_4$ on stationary Bratteli diagrams such that $S(\mu_3) = S(\mu_4)$ but $\mu_3$ and $\mu_4$ are not homeomorphic. Therefore, the invariant $S(\mu)$ is not a complete invariant for homeomorphism equivalence of measures.

\unitlength = 0.4cm
\begin{center}
\begin{tabular*}{0.99\textwidth}%
     {@{\extracolsep{\fill}}cc}

\begin{graph}(14,13)
\graphnodesize{0.4}
\roundnode{V0}(7,12)
\roundnode{V10}(1,9)
\roundnode{V11}(5,9)
\roundnode{V12}(9,9)
\roundnode{V13}(13,9)
\edge{V10}{V0}
\edge{V11}{V0}
\edge{V12}{V0}
\edge{V13}{V0}
\roundnode{V20}(1,4.5)
\roundnode{V21}(5,4.5)
\roundnode{V22}(9,4.5)
\roundnode{V23}(13,4.5)
\bow{V20}{V10}{0.04}
\bow{V20}{V10}{0.12}
\bow{V20}{V10}{-0.12}
\bow{V20}{V10}{-0.04}
\edge{V21}{V10}
\bow{V21}{V11}{0.06}
\bow{V21}{V11}{-0.06}
\edge{V22}{V11}
\bow{V22}{V12}{0.10}
\edge{V22}{V12}
\bow{V22}{V12}{-0.10}
\edge{V22}{V13}
\bow{V23}{V13}{0.10}
\edge{V23}{V13}
\bow{V23}{V13}{-0.10}
\edge{V23}{V11}
\edge{V23}{V12}
\roundnode{V30}(1,0.5)
\roundnode{V31}(5,0.5)
\roundnode{V32}(9,0.5)
\roundnode{V33}(13,0.5)
\bow{V30}{V20}{0.04}
\bow{V30}{V20}{0.12}
\bow{V30}{V20}{-0.12}
\bow{V30}{V20}{-0.04}
\edge{V31}{V20}
\bow{V31}{V21}{0.06}
\bow{V31}{V21}{-0.06}
\edge{V32}{V21}
\bow{V32}{V22}{0.10}
\edge{V32}{V22}
\bow{V32}{V22}{-0.10}
\edge{V32}{V23}
\bow{V33}{V23}{0.10}
\edge{V33}{V23}
\bow{V33}{V23}{-0.10}
\edge{V33}{V21}
\edge{V33}{V22}
\freetext(7,-0.9){$.\ .\ .\ .\ .\ .\ .\ .\ .\ .\ .\ .\ .\ .\ .\ .\ .\ .\ .$}%
 \freetext(5,-2.5){{$B_3$: $\mu_3$ is good}}
\end{graph}

            &
\unitlength = 0.4cm
\begin{graph}(14,13)
\graphnodesize{0.4}
\roundnode{V0}(7,12)
\roundnode{V10}(1,9)
\roundnode{V11}(5,9)
\roundnode{V12}(9,9)
\roundnode{V13}(13,9)
\edge{V10}{V0}
\edge{V11}{V0}
\edge{V12}{V0}
\edge{V13}{V0}
\roundnode{V20}(1,4.5)
\roundnode{V21}(5,4.5)
\roundnode{V22}(9,4.5)
\roundnode{V23}(13,4.5)
\bow{V20}{V10}{0.09}
\bow{V20}{V10}{0.18}
\edge{V20}{V10}
\bow{V20}{V10}{-0.18}
\bow{V20}{V10}{-0.09}
\edge{V21}{V10}
\bow{V21}{V11}{0.06}
\bow{V21}{V11}{-0.06}
\edge{V22}{V11}
\bow{V22}{V12}{0.04}
\bow{V22}{V12}{0.12}
\bow{V22}{V12}{-0.12}
\bow{V22}{V12}{-0.04}
\edge{V22}{V13}
\bow{V23}{V13}{0.04}
\bow{V23}{V13}{0.12}
\bow{V23}{V13}{-0.12}
\bow{V23}{V13}{-0.04}
\edge{V23}{V11}
\edge{V23}{V12}
\roundnode{V30}(1,0.5)
\roundnode{V31}(5,0.5)
\roundnode{V32}(9,0.5)
\roundnode{V33}(13,0.5)
\bow{V30}{V20}{0.09}
\bow{V30}{V20}{0.18}
\edge{V30}{V20}
\bow{V30}{V20}{-0.18}
\bow{V30}{V20}{-0.09}
\edge{V31}{V20}
\bow{V31}{V21}{0.06}
\bow{V31}{V21}{-0.06}
\edge{V32}{V21}
\bow{V32}{V22}{0.04}
\bow{V32}{V22}{0.12}
\bow{V32}{V22}{-0.12}
\bow{V32}{V22}{-0.04}
\edge{V32}{V23}
\bow{V33}{V23}{0.04}
\bow{V33}{V23}{0.12}
\bow{V33}{V23}{-0.12}
\bow{V33}{V23}{-0.04}
\edge{V33}{V21}
\edge{V33}{V22}
\freetext(7,-0.9){$.\ .\ .\ .\ .\ .\ .\ .\ .\ .\ .\ .\ .\ .\ .\ .\ .\ .\ .$}%
 \freetext(5,-2.5){{$B_4$: $\mu_4$ is not good}}
\end{graph}

\vspace{1.5cm}
\end{tabular*}
\end{center}

\noindent\textbf{Example 2.} We show that there are non-ergodic infinite $\mathcal{R}$-invariant measures on stationary Bratteli diagrams that belong to $M_\infty(X) \setminus M_\infty^{0}(X)$. Let $B$ be a stationary Bratteli diagram with the matrix $A = F^T$ such that
$$
A =
\begin{pmatrix}
A_1 & X_{12} & X_{13}\\
0 & A_2 & X_{23}\\
0 & 0 & A_3 \\
\end{pmatrix},
$$
where the square non-zero matrices $\{A_i\}_{i=1}^3$ are strictly positive, $X_{23}$ is non-zero and at least one of the matrices $X_{12}$, $X_{23}$ is non-zero.
Then the diagram $B$ has three subdiagrams defined by matrices $\{A\}_{i=1}^3$ and three ergodic $\mathcal{R}$-invariant measures.
Denote by $\alpha_i$ the class of vertices corresponding to $A_i$ for $i = 1,2,3$. Denote by $\mu_1$, $\mu_2$, $\mu_3$ the corresponding ergodic $\mathcal{R}$-invariant measures on $B$. Recall that each measure $\mu_i$ is full on its support $\overline{X_{\alpha_i}}$ for $i = 1,2,3$. Note that $\overline{X_{\alpha_1}} = X_{\alpha_1}$ and $\overline{X_{\alpha_3}} = X_B$. The measure $\mu_1$ is finite while $\mu_i \in M^0_\infty(\overline{X_{\alpha_i}})$ for $i = 2, 3$.

Let the spectral radii of the matrices $A_i$ satisfy the following inequalities: $sp(A_1) > sp(A_2) > sp(A_3)$. Let $\nu_{12} = \mu_1 + \mu_2$ and $\nu_{123} = \mu_1 + \mu_2 + \mu_3$.
Then $\nu_{12} \in M_\infty(\overline{X_{\alpha_2}})$ and $\nu_{123} \in M_\infty(X_B)$.
We have $\mathfrak{M}_{\nu_{12}} = X_{\alpha_1}$ hence $\nu_{12}(\mathfrak{M}_{\nu_{12}}) = \mu_1 (X_{\alpha_1}) \in (0, \infty)$. Note that the values of $\mu_2$ and $\nu_{12}$ on the clopen sets are the same.
Similarly, $\mathfrak{M}_{\nu_{123}} = \overline{X_{\alpha_2}}$ and $\nu_{123}(\mathfrak{M}_{\nu_{123}}) = \mu_1 (X_{\alpha_1}) + \mu_2 (\overline{X_{\alpha_2}}) = \infty$.

\medskip
The next proposition can be proved the same way as Proposition 5.3 in~\cite{S.B.O.K.}.
\begin{prop}
Let $\mu$ be any infinite ergodic $\mathcal R$-invariant measure on a stationary Bratteli diagram. Then there exists a good infinite ergodic $\mathcal R$-invariant measure $\nu$ on a stationary Bratteli diagram such that $S(\mu) = S(\nu)$ and $\mathfrak{M}$ is homeomorphic to $\mathfrak{N}$.
\end{prop}

\Large \noindent \textbf{Acknowledgement}

\medskip
\normalsize
I would like to thank my advisor Sergey Bezuglyi for introducing this subject to me, for numerous helpful discussions, for providing me with the idea of compactification of a locally compact space (Theorem~\ref{goodSmu}, Corollary~\ref{invarhomeo}) to obtain a good non-defective measure and for reading the preliminary versions of this paper.

\end{document}